\documentclass[preprint,12pt]{elsarticle}
\usepackage{geometry}
\usepackage{amssymb, amsmath}
\usepackage{amsthm}
\usepackage[utf8]{inputenc}

\usepackage{hyperref}
\usepackage{enumitem}
\usepackage{tabularx}
\hypersetup{
  colorlinks=true,
}
\hypersetup{
    linkcolor=magenta,
    filecolor=magenta,      
    urlcolor=cyan,
    citecolor=magenta,
    }

\theoremstyle{remark}

\theoremstyle{definition}
\newtheorem{theorem}{Theorem}

\journal{xxxxxxxxxxx}

\newcommand{\cinf}{$\mathcal{C}^{\infty}$}

\DeclareMathOperator{\csch}{csch}
\DeclareMathOperator{\jacobisn}{sn}
\DeclareMathOperator{\ellipticF}{F}
\DeclareMathOperator{\arctanh}{arctanh}
\DeclareMathOperator{\arccot}{arccot}
\DeclareMathOperator{\sech}{sech}

\usepackage{todonotes}
\usepackage{xcolor}

\usepackage{marginnote}
\usepackage{tikz}
\begin{document}
\hypersetup{
    linkcolor=magenta,
    filecolor=magenta,      
    urlcolor=cyan,
    citecolor=magenta,
    }
\begin{frontmatter}



\title{Classification of traveling wave solutions of the modified Zakharov--Kuznetsov equation}


\author[pan]{Antonio J. Pan-Collantes}

\author[pan,conchi]{C. Muriel}

\author[pan]{A. Ruiz}

\affiliation[pan]{
organization={Departamento de Matematicas,
Universidad de Cadiz - UCA},%Department and Organization
            addressline={Facultad de Ciencias, Campus Universitario de Puerto Real s/n}, 
            city={Puerto Real},
            postcode={11510}, 
            state={Cadiz},
            country={Spain}}

\affiliation[conchi]{
organization={International Society of Nonlinear Mathematical Physics (ISNMP)},%Department and Organization
            addressline={Auf der Hardt 27}, 
            city={Bad Ems},
            postcode={56130}, 
            country={Germany}}

\begin{abstract}
The \cinf-structure-based method of integration of distributions of vector fields is used to classify all the traveling wave solutions of the modified Zakharov--Kuznetsov equation. This work unifies and generalizes the particular results obtained in the recent literature by using specific ansatz-based methods.
\end{abstract}


\end{frontmatter}


\section{Introduction}
%%%%%%%%%%%%%%%%%%%%%%%%%%%%%%%%%%%%%%%%%%%%%%

Nonlinear partial differential equations play a fundamental role in the study of a wide range of physical phenomena. They can be used to model the behavior of a system with respect to time (evolution equation), taking into account nonlinear effects such as dispersion, dissipation, and non-local interactions. 

A well-known example of nonlinear evolution equation is the Zakharov--Kuznetsov equation, used to study the behavior of weakly nonlinear ion-acoustic waves within a plasma comprising hot electrons and cold ions in the presence of an unvarying magnetic field oriented along the $x$-axis. It was first presented by Zakharov and Kuznetsov \cite{ZKoriginal}, who derived the following equation for $u=u(x,y,t)$:
\begin{equation}\label{c6EqZK_orig}
	u_t+u u_x+\nabla^2 u_x=0,
\end{equation}
where $\nabla^2=\frac{\partial^2}{\partial x^2} + \frac{\partial^2}{\partial y^2}$ is the Laplacian operator.

This equation has been found to have plenty of other applications ranging from Rossby waves in a rotating atmosphere to long waves existing on a thin liquid film, and even to the isolated vortex of drift waves in a three-dimensional plasma. Thus, it is widely applied in various fields of science, including optical fiber, geochemistry, and solid-state physics (see \cite{lagrangianZK} and references therein).

In this work, we discuss the modified Zakharov--Kuznetsov equation (mZK), as referenced in \cite{lagrangianZK,mZKOsman,gZKsymmetry,gZKvarcoef}:
\begin{equation}\label{c6EqGenZK}
	u_t+A u u_x+B u^2 u_x+M u_{xxx}+Nu_{xyy}=0,
\end{equation}
where $A,B,M,N$ are real constants. In particular, we will focus on the analysis of the traveling wave solutions \cite{olver86,blumanlibro} of this equation, which represent solutions that propagate in time while maintaining their shape. These solutions often arise in physical and biological systems where waves play a fundamental role, such as in fluid mechanics, chemical reactions, and neural networks. The study of traveling wave solutions can provide insight into the stability and dynamics of the system under consideration, and can also be used to predict the onset of important phenomena such as pattern formation and wave breaking.

In the recent literature there has been an intense research activity regarding the determination of exact solutions to different types of Zakharov--Kuznetsov equations. For instance, in 
\cite{mZKOsman} M. A. Iqbal {\it et al.} obtained particular traveling wave  solutions to the mZK equation \eqref{c6EqGenZK} by using the $(G'/G,1/G)$-expansion method; in \cite{lagrangianZK} C. M. Khalique and O. D. Adeyemo studied a $(3+1)$-dimensional generalized Zakharov--Kuznetsov equation (gZK) with dual power law nonlinearity and obtained particular exact solutions via the Lie group method as well as the direct integration method and the extended Jacobi elliptic function expansion approaches; in  \cite{gZKsymmetry}  Y. Gu and J. Qi also computed particular solutions to the gZK equation by using the $e^{-\phi(z)}$-expansion method; in \cite{gZKvarcoef} S. K. Mohanty {\it et al.} dealt with the gZK equation in which the coefficients depend on the time and particular solutions by means of the $G/G'$-expansion method were obtained; particular solutions to equation \eqref{c6EqGenZK} for the case $A=0$ and $B=M=N=1$ were obtained in  \cite{TASCAN20091810} by F. Tascan {\it et al.} using the first integral method  and in \cite{ESLAMI2014221} by M. Eslami {\it et al.}  by means of the homogeneous balance method, whereas other particular solutions to equation \eqref{c6EqGenZK} for $B=0$ and $M=N$ were computed by Y. Yu {\it el al.} in \cite{yu}. Notwithstanding these significant contributions, a complete classification of all the traveling wave solutions to equation \eqref{c6EqGenZK}  has not been performed yet in the literature.

Therefore, the aim of this paper is to provide a complete classification of all the traveling wave type solutions to the mZK equation \eqref{c6EqGenZK}. For this purpose, the \cinf-structure-based method, recently presented in \cite{pancinf-sym,pancinf-struct,pan23integration}, is applied. This integration method allows to determine the integral manifolds of an involutive distribution of vector fields and, in particular, to find solutions to $m$th-order ordinary differential equations (ODEs) by solving $m$ completely integrable Pfaffian equations.

The paper is organized as follows: in Section \ref{travelingwaves} the traveling wave reduction is applied to equation \eqref{c6EqGenZK}, leading to a third-order ODE. After that, a \cinf-structure for the (trivially) involutive distribution spanned by the vector field 
associated to the ODE is computed. The application of the \cinf-structure-based integration method permits to obtain all the traveling wave solutions of the mZK equation \eqref{c6EqGenZK}, which are expressed in implicit form and in terms of a primitive involving a square root of a rational function. In Section \ref{classification}, a classification of the traveling wave solutions, which can be expressed in explicit form, is carried out, according to the different roots and their multiplicities of the denominator of the corresponding rational function. In Section \ref{section4}, it is described how particular solutions obtained in the previous literature can fit into the complete classification performed in this work. Finally, in Section \ref{section5} some particular examples of the classified explicit solutions are showed. By selecting adequately the parameters of the equation and the integration constants, different types of traveling waves with physical interest can be obtained, including bright solitons, kink and periodic solutions.

%%%%%%%%%%%%%%%%%%%%%%%%%%%%%%
\section{Traveling wave solutions of the mZK equation}\label{travelingwaves}
%%%%%%%%%%%%%%%%%%%%%%%%%%%

% maple 098 and mathematica 011

In order to study the traveling wave solutions of the modified Zakharov--Kuznetsov equation \eqref{c6EqGenZK} we apply the transformation
\begin{equation}\label{c6EqTrav2D}
	\begin{array}{l}
		r=x+y-c t, \quad (c\neq 0),\\	
		v(r)=u(x,y,t),	
	\end{array}
\end{equation}
 to equation \eqref{c6EqGenZK}, and we obtain the third-order ordinary differential equation (ODE)
\begin{equation}\label{ZKode}
-cv_1+Avv_1+Bv^2 v_1+(M+N)v_3=0,
\end{equation}
where $v=v(r)$ and $v_i$ denotes the $i$th-order derivative of $v$ with respect to $r$, for $1\leq i \leq 3$.

We assume $M+N\neq 0$ because the case $M+N=0$ corresponds to 
\begin{equation}\label{ZKodefactoriz}
	(-c+Av+Bv^2) v_1=0,
	\end{equation}
which only admits constant solutions $v(r)=K$ for $K\in \mathbb R$.

Equation \eqref{c6EqGenZK} may be reduced to a second-order ODE, either through direct integration or by using the admitted Lie point symmetry $\partial_r$. In both cases, if we find a Lie point symmetry of the reduced equation (which would be a Type II hidden symmetry \cite{abraham2002hidden} of the original equation),  the integration process can be continued by reducing to a first-order ODE. Eventually, once the latter is integrated, solutions of the original third-order ODE \eqref{ZKode} could be recovered by inverting several changes of variables and performing quadratures \cite{olver86,ovsiannikovlibro,stephani}.

Alternatively, the integration of \eqref{ZKode} can be achieved by using the recent \cinf-structure-based method \cite{pancinf-sym,pancinf-struct,pan23integration}. 
This integration process is carried out, as in the case of solvable structures \cite{basarab,hartl1994solvable,sherring1992geometric}, by solving a sequence of Pfaffian equations written in the same variables as the original equation \eqref{ZKode}. In this way, we avoid inverting the changes of variables and performing the quadratures to recover the solution that appears in the classical Lie method.

In order to find exact solutions to \eqref{ZKode} by using the \cinf-structure-based method, we consider the jet space $J^{2}(\mathbb R, \mathbb R)$, with coordinates $(r,v,v_1,v_2)$. The vector field
$$
Z=\partial_r+v_1 \partial_v+v_2\partial_{v_1}-\dfrac{(Bv^2+Av-c)v_1}{M+N}\partial_{v_2},
$$
defined on $J^{2}(\mathbb R, \mathbb R)$, represents equation \eqref{ZKode}, in the sense that its integral curves are in correspondence with the prolongations of the solutions of \eqref{ZKode}. The procedure consists of the following steps (the details can be found in \cite{pancinf-sym, pancinf-struct,pan23integration}):
\begin{enumerate}
    \item Determination of a \cinf-structure $\langle X_1, X_2,X_3\rangle$: In order to apply the procedure, we need to identify a \cinf-structure for the involutive distribution generated by $Z$, i.e., an ordered set of vector fields $\langle X_1, X_2,X_3\rangle$ such that $\{Z,X_1,X_2,X_3\}$ are pointwise linearly independent and the distribution generated by $\{Z,X_1\}$ and the distribution generated by $\{Z,X_1,X_2\}$ are both involutive.
    
    The vector field $\partial_r$ is a Lie point symmetry of \eqref{ZKode}, so it can be used as the first vector field  in our \cinf-structure: $X_1:=\partial_r$.

To find the second vector field $X_2$ in our \cinf-structure, we can consider the following {\it ansatz}:
$
X_2=\partial_{v_1}+\eta \partial_{v_2},
$
where $\eta=\eta(r,v,v_1)$ is a smooth function to be determined by imposing that  $[X_2,Z]$ and $[X_2,X_1]$ belong to the distribution generated by $\{Z,X_1,X_2\}.$ The following system of determining equations for $\eta=\eta(r,v,v_1)$ is obtained:
\begin{equation}
	\begin{aligned}
	v_1\eta_r +v_1^2 \eta_v+v_2 v_1 \eta_{v_1}-v_2 \eta+v_1\eta^2=0,\\
	-\eta_r=0.		
	\end{aligned}
\end{equation}

The particular solution $\eta=\dfrac{v_1}{v}$ can be easily found.

The third element $X_3$ of the \cinf-structure can be any vector field such that the set $\{Z,X_1,X_2,X_3\}$ is pointwise linearly independent; we can choose, for simplicity, $X_3=\partial_{v_2}$.

It can be checked that  the vector fields \begin{equation}\label{c6EqCanZK2}
X_1=\partial_r,\quad
		X_2=\partial_{v_1}+\dfrac{v_1}{v}\partial_{v_2}, \quad
		X_3=\partial_{v_2}
\end{equation} satisfy the following commutation relations 
$$
\begin{aligned}
&[X_1, Z] = [X_1, X_2] = [X_1, X_3] = [X_2, X_3] = 0, \\
&[X_2, Z] = \frac{1}{v_1}(Z - X_1) - \frac{v v_2 - v_1^2}{v v_1} X_2, \\
&[X_3, Z] = X_2 - \frac{v_1}{v} X_3.
\end{aligned}
$$
Therefore $\langle X_1,X_2,X_3\rangle$ is a \cinf-structure for equation \eqref{ZKode}.
\item Calculation of the 1-forms that define the sequence of Pfaffian equations: Let $\boldsymbol{\Omega} = dr \wedge dv \wedge d v_1 \wedge dv_2$  denote the volume form in $J^2(\mathbb{R},\mathbb{R})$ and $\lrcorner$ the interior product (or contraction).  Then, following the procedure outlined in \cite[Section 3]{pan23integration}, we compute the 1-forms:
\begin{equation}\label{c6EqZKforms}
	\begin{aligned}
		\omega_1&={X_{3}\,\lrcorner\,X_{2}\,\lrcorner\,Z\,\lrcorner\,\boldsymbol{\Omega}}%\\ &
  =-v_1 dr+dv,\quad \\[0.4cm]
		\omega_2&={X_{3}\,\lrcorner\,X_{1}\,\lrcorner\,Z\,\lrcorner\,\boldsymbol{\Omega}}%\\&
		=-v_2 dv + v_1 dv_1, \\[0.4cm]
		\omega_3&={X_{2}\,\lrcorner\,X_{1}\,\lrcorner\,Z\,\lrcorner\,\boldsymbol{\Omega}}=\\&
		=-\dfrac{\left(B v^3+A v^2+(M+N)v_2-cv\right)v_1 }{(M+N) v}dv+\dfrac{v_1^2}{v}dv_1- v_1 dv_2.
	\end{aligned}	
	\end{equation}

\item Pfaffian equations and successive integration:  
\begin{itemize}
    \item  The Pfaffian equation $\omega_3\equiv 0$ is completely integrable. A corresponding first integral $I_3=I_3(r,v,v_1,v_2)$ satisfies $dI_3\wedge \omega_3=0,$ which yields to a homogeneous linear system of partial differential equations for  $I_3$:
\begin{equation}
\begin{aligned}
    &(I_3)_r = 0, \\[10pt]
     &(B v^3+Av^2-cv+ (M+N) v_2)(I_3)_{v_1}   + (M+N) v_1(I_3)_v  = 0, \\[10pt]
     &v_1(I_3)_{v_2}  +v (I_3)_{v_1}   = 0, \\[10pt]
    &\frac{ v_1 (B v^3 + A v^2 + (M+N) v_2  - c v)}{v (M + N)}(I_3)_{v_2} - v_1 (I_3)_v = 0.
\end{aligned}    
\end{equation}
It can be checked that a particular solution for that system becomes
\begin{equation}\label{c6EqZKI3}
I_3=3Bv^4+4Av^3-6cv^2+6(M+N)(2vv_2-v_1^2).
\end{equation}

\item For $C_3\in \mathbb{R}$ we consider the corresponding level set of $I_3:$   
$$
\Sigma_{(C_3)}=\{(t,v,v_1,v_2)\in J^2(\mathbb{R},\mathbb{R}):
I_3(r,v,v_1,v_2)=C_3\},
$$
which can be locally parametrized by  
$$
\iota_3(r,v,v_1) =\left(r,v,v_1,-\dfrac{3Bv^4+4Av^3-6(M+N)v_1^2-6cv^2-C_3}{12\left(M+N\right)v}\right).
$$
The restrictions of the forms \eqref{c6EqZKforms} to $\Sigma_{(C_3)}$, i.e., $\omega_i|_{\Sigma_{(C_3)}}:=\iota_3^*(\omega_i)$ for $i=1,2,3$, become:
\begin{equation}
\begin{aligned}\label{c6EqZKforms2}
	\omega_1|_{\Sigma_{(C_3)}}&=-v_1 dr+dv,\\
	\omega_2|_{\Sigma_{(C_3)}}&=\dfrac{3Bv^4+4Av^3-6(M+N)v_1^2-6cv^2-C_3}{12\left(M+N\right)v}dv+v_1dv_1,\\
	\omega_3|_{\Sigma_{(C_3)}}&=0.
\end{aligned}
\end{equation}

According to \cite[Theorem 3.5]{pancinf-sym}, the Pfaffian equation  $\omega_2|_{\Sigma_{(C_3)}}\equiv 0$ is completely integrable. As before, a corresponding primitive $I_2=I_2(r,v,v_1;C_3)$ satisfies $dI_2\wedge \omega_2|_{\Sigma_{(C_3)}}=0.$ It can be checked that 
a particular solution for the corresponding system of determining equations is the smooth function
\begin{equation}\label{c6EqZKI2}
	I_2=\dfrac{Bv^4+2Av^3+6(M+N)v_1^2-6cv^2+ C_3}{v}, \, v\neq 0.
\end{equation}
\item For $C_2\in \mathbb{R}$ let ${\Sigma_{(C_2,C_3)}}$ denote the level set defined by $I_2(r,v,v_1;C_3)=C_2.$ 
A local parametrization of ${\Sigma_{(C_2,C_3)}}$ is  given by
\begin{equation}\label{iota2ZK}
	\iota_2(r,v) = \left(r,v,\mp\sqrt{\dfrac{-Bv^4-2Av^3+6cv^2+C_2v-C_3}{6(M+N)}}\right).
\end{equation}
The restrictions of the 1-forms \eqref{c6EqZKforms2}  to ${\Sigma_{(C_2,C_3)}}$ via $\iota_2^*$ become
\begin{equation}\label{unoformafinal}
	\begin{aligned}
\omega_1|_{\Sigma_{(C_2,C_3)}}&=\pm\sqrt{\dfrac{-Bv^4-2Av^3+6cv^2+C_2v-C_3}{6(M+N)}} dr+dv,\\
		\omega_2|_{\Sigma_{(C_2,C_3)}}&=0,\\
		\omega_3|_{\Sigma_{(C_2,C_3)}}&=0.
	\end{aligned}	
\end{equation}

As previously, the Pfaffian equation $\omega_1|_{\Sigma_{(C_2,C_3)}}\equiv 0$ is completely integrable. Furthermore, this Pfaffian can be solved by quadrature:   %Indeed, $\sqrt{\frac{6(M+N)}{-Bv^4-2Av^3+6cv^2+C_2v-C_3}}$ is an integrating factor. 
a corresponding first integral  $I_1=I_1(r,v;C_2,C_3)$   is formally given by
$$
I_1=r\pm H(v;C_2,C_3),
$$
where $H(v;C_2,C_3)$ is any primitive of the function
\begin{equation}\label{integral}
	h(v;C_2,C_3)=\sqrt{\dfrac{-6(M+N)}{B v^4+2Av^3-6cv^2-C_2v+C_3}}.
\end{equation} 
\end{itemize}
\end{enumerate}

Therefore, the preceding discussion yields the following result:

\begin{theorem}
All traveling wave solutions of the mZK equation \eqref{c6EqGenZK} are implicitly described by
\begin{equation}\label{implicitTodas}
	x + y - ct \pm H(u(x,y,t);C_2,C_3) = C_1,
\end{equation}
where $C_1,C_2,C_3 \in \mathbb{R}$ and $H(v;C_2,C_3)$ is a primitive of the function
\begin{equation}\label{h}
h(v;C_2,C_3) = \sqrt{\dfrac{-6(M+N)}{B v^4 + 2 A v^3 - 6 c v^2 - C_2 v + C_3}}.
\end{equation}
\end{theorem}

In the following section we  analyze the implicit equation \eqref{implicitTodas}, in order to obtain a complete classification of the traveling wave solutions of the mZK equation in closed form.

\section{Classification of all traveling wave solutions in explicit form}
\label{classification}
%%%%%%%%%%%%%%%%%%%%%%%%%%%%%%%%%%%%%%
%%%%%%%%%%%%%%%%%%%%%%%%%%%%%%%%%%%%%%%%%%%%%5
In this section we first investigate how to obtain explicit expressions for the smooth function $H$ that appears  in \eqref{implicitTodas}. Such function $H=H(r,v;C_2,C_3)$ is a primitive of the function $h(v;C_2,C_3)$ given in \eqref{h}. For this purpose, we will consider all the possible cases depending on the roots of the polynomial
\begin{equation}\label{polinomioP}
	P(z)=B z^4+2Az^3-6cz^2-C_2z+C_3
\end{equation}
and their multiplicities.

Throughout this section, we will assume that the domains (when not explicitly specified) are properly constrained to ensure that functions and expressions are well defined.

\subsection{\textbf{Case I}: \texorpdfstring{$A=B=0$}{}}
%%%%%%%%%%%%%%%%%%%%%%%%%%%%%%%%%%%%%%%%%%
%ver mathematica 021

    When $A=B=0$ the polynomial $P(z)$ has degree two and the function $h$ in \eqref{h} reduces to
\begin{equation}\label{hI}
    h(v)=\sqrt{\dfrac{6(M+N)}{6cv^2+C_2v-C_3}}.
\end{equation}
	
	We distinguish two cases:
	\begin{enumerate}
		\item If $\frac{M+N}{c}>0$, then necessarily $v$ is such that $\frac{c}{6cv^2+C_2v-C_3}>0$, and a primitive of \eqref{hI} is given by
		$$
H(v)=-\sqrt{\dfrac{M+N}{c}}\ln \left(\sqrt{\frac{6cv^2+C_2v-C_3}{9c}}-2v-\frac{C_2}{6 c}\right).
		$$
		The corresponding solutions to the modified ZK equation \eqref{c6EqGenZK}, obtained from \eqref{implicitTodas}, become
\begin{equation}\label{solg21}
	u(x,y,t) =-\frac{C_2}{12 c}-\frac{1}{4}e^{\mp\sqrt{\frac{c}{M+N}}(C_1-r)}-\frac{C_2^2+24cC_3}{144c^2}e^{\pm\sqrt{\frac{c}{M+N}}(C_1-r)},
\end{equation}
where $r=x+y-ct$. 

	\item If $\frac{M+N}{c}<0$, then necessarily $v$ is such that $\frac{c}{6cv^2+C_2v-C_3}<0$, and a primitive \eqref{hI} is given by
	$$
H(v)=2\sqrt{-\frac{M+N}{c}}\arctan\left(\frac{v}{-\sqrt{\frac{C_3}{6c}}+\sqrt{\frac{-6cv^2-C_2v+C_3}{6c}}}\right).
	$$
	The corresponding solutions to the modified ZK equation \eqref{c6EqGenZK}, obtained from \eqref{implicitTodas}, is 
\begin{equation}\label{solg22}
	u(x,y,t) =\frac{ \mp\sqrt{\frac{24 C_3}{c}}\tan\left(\xi\right)-\frac{C_2}{c}\tan^2\left(\xi\right)    }{6 \left(1+\tan^2\left(\xi\right)\right)}      ,
\end{equation}
where $\xi=\sqrt{\frac{-c}{4 (M+N)}} (C_1-x-y+ct)$.
	\end{enumerate}
In expressions \eqref{solg21} and \eqref{solg22} the parameters $c,C_1,C_2,C_3$ are arbitrary.

\subsection{\textbf{Case II}: \texorpdfstring{$A\neq 0$ and $B=0$}{}}
%%%%%%%%%%%%%%%%%%%%%%%%%%%%%%%%%%%%%%%%%%%%%%%%%%
%maple 098 y mathematica 022
In this case, the polynomial $P(z)$ is cubic, and we distinguish three cases:
\begin{itemize}
	\item \textbf{Case II.1: $P(z)$ has a triple root}. This situation arises if and only if the coefficients satisfy: 
	$$
	C_2=\frac{-6c^2}{A}, \quad	C_3=\frac{-2 c^3}{A^2}.
	$$ 
	The triple root is $\varphi_1=\frac{c}{A}$, and the function $h$ takes the form
	$$
	h(v)=  \sqrt{\dfrac{3A^2(M+N)}{  (c-Av)^3}}.
	$$

	A primitive is given by
	$$
	H(v)=\frac{2}{A}\sqrt{\dfrac{3A^2(M+N)}{  (c-Av)^3}} (c-Av),
	$$
	which yields to the following two-parameter family of solutions for the mZK equation:
\begin{equation}\label{solg3triple}
	u(x,y,t)=- \frac{12(M+N)}{A(C_1-r)^2}+\frac{c}{A},
\end{equation}
	where $r=x+y-ct$. Note that the parameters $c$ and $C_1$ are free.

	\item \textbf{Case II.2: $P(z)$ has a double root $\varphi$}. In this case we can write
    $$
	C_2=6 A \varphi^2-12 c \varphi,\quad C_3=4A\varphi^3-6c\varphi^2,
	$$ where  $\varphi$ is the double root and the single root is $\tilde{\varphi}=\frac{3 c}{A}-2 \varphi.$ In this way, the function
 $h$ in \eqref{h} can be expressed as
    \begin{equation}\label{g3raizdoble}
		h(v;\varphi) = \sqrt{\dfrac{3(M+N)}{  (v-\varphi)^2(3 c-Av-2A\varphi)}}.	
	\end{equation}

To find a primitive for \eqref{g3raizdoble} we distinguish the following subcases:
\begin{enumerate}[label=(\alph*)]
	\item If $M+N>0$ then $v$ is necessarily such that $3 c-Av-2A\varphi>0$. %$\tilde{\varphi}-v>0$
	When $c-A\varphi>0$ a primitive of  \eqref{g3raizdoble}  is given by
	\begin{equation}\label{primitivaatanh}
		H(v;\varphi)=\pm 2 \sqrt{\frac{M+N}{c-A\varphi}} \arctanh\left(\sqrt{\frac{3c-2A\varphi-Av}{3c-3A\varphi}}\right)
	\end{equation}
	 %$\varphi<\tilde{\varphi}$ 
	The choice of the $+$ or $-$ sign corresponds to whether $v < \varphi$ or $v > \varphi$, respectively.

	If $c-A\varphi<0$ %$\varphi>\tilde{\varphi}$  
	then necessarily is $v<\varphi$, and a primitive of  \eqref{g3raizdoble} is given by 
	$$
H(v;\varphi)=- 2 \sqrt{\frac{M+N}{A\varphi-c}} \arctan\left(\sqrt{\frac{3c-2A\varphi-Av}{3A\varphi-3c}}\right).
	$$

The corresponding families of solutions to equation \eqref{c6EqGenZK} are given, respectively, by
\begin{equation}\label{solg3caso2a1}
u(x,y,t)=\varphi+\frac{6(c-A\varphi)}{A+A\cosh\left(\sqrt{\frac{c-A\varphi}{M+N}}(C_1-r)\right)},
\end{equation}
and
\begin{equation}\label{solg3caso2a2}
	u(x,y,t)=\varphi+\frac{6(A\varphi-c)}{A+A\cos\left(\sqrt{\frac{A\varphi-c}{M+N}}(C_1-r)\right)},
	\end{equation}
where $r=x+y-ct$. In \eqref{solg3caso2a1} $c, \varphi$ and $C_1$ are three arbitrary parameters.

\item If $M+N<0$ then $v$ is necessarily such that $3 c-Av-2A\varphi<0$. %$\tilde{\varphi}-v<0$
If additionally $c-A\varphi>0$, then a primitive of  \eqref{g3raizdoble}  is given by
\begin{equation}\label{primitArcoCot}
	H(v;\varphi)=2\sqrt{\frac{-(M+N)}{c-A\varphi}} \arccot\left(\sqrt{\frac{3c-3A\varphi}{Av+2A\varphi-3c}}\right).
\end{equation}

When $c-A\varphi<0$, a primitive of  \eqref{g3raizdoble}  is given by the same expression \eqref{primitivaatanh}.

The family of solutions for equation \eqref{c6EqGenZK} obtained from \eqref{primitArcoCot} is
\begin{equation}\label{solg3caso2b}
	u(x,y,t)=\varphi+\frac{3c-3A\varphi}{A\cos^2\left(\sqrt{\frac{A\varphi-c}{4(M+N)}}(C_1-r)\right)},
\end{equation}
where $r=x+y-ct$ and $c,C_1,\varphi$ are arbitrary constants.

\end{enumerate}

	\item \textbf{Case II.3: $P(z)$ has three simple roots $\varphi_1,\varphi_2$ and $\varphi_3$}.
	This case arises exactly when
	$$
	C_2=2A\varphi_1^2+2A\varphi_1 \varphi_2+2A\varphi_2^2-6c\varphi_1-6c\varphi_2
	$$
	and
	$$
C_3=2A\varphi_1^2\varphi_2+2A\varphi_1\varphi_2^2-6c\varphi_1\varphi_2,
	$$
for $\varphi_1,\varphi_2 \in \mathbb R$, with $\varphi_1\neq \varphi_2,\varphi_1\neq \frac{1}{2}(\frac{3c}{A}-\varphi_2)$ and $\varphi_2\neq \frac{1}{2}(\frac{3c}{A}-\varphi_1)$. The roots of the polynomial $P(z)$ are $\varphi_1,\varphi_2$ and $\varphi_3= \frac{3 c}{A}-\varphi_1-\varphi_2$. 

	In this case, we have to look for a primitive of
	$$
h(v;\varphi_1,\varphi_2)=\sqrt{\dfrac{-3(M+N)}{A \left(v-\varphi_1\right)\left(v-\varphi_2\right)\left(v+  \varphi_1+\varphi_2-\frac{3 c}{A}\right)}}.
	$$

	We can assume, without loss of generality, that $\varphi_1<\varphi_2<\varphi_3$. We  distinguish the following subcases:
	\begin{enumerate}[label=(\alph*)]
		\item If $\frac{M+N}{A}<0$ then necessarily $v>\varphi_3$ or $\varphi_1<v<\varphi_2$. By performing, respectively, the change of variables $\tilde{v}=\sqrt{\frac{\varphi_2-\varphi_1}{v-\varphi_1}}$ or $\tilde{v}=\sqrt{\frac{v-\varphi_1}{\varphi_2-\varphi_1}}$, the following primitives are obtained
		$$
H(v;\varphi_1,\varphi_2)=-\sqrt{\frac{-12(M+N)}{A(\varphi_2-\varphi_1)}} \ellipticF\left(\sqrt{\frac{\varphi_2-\varphi_1}{v-\varphi_1}},\frac{\varphi_3-\varphi_1}{\varphi_2-\varphi_1}\right),\text{ for } v>\varphi_3,
		$$
		and
		$$
H(v;\varphi_1,\varphi_2)=\sqrt{\frac{-12(M+N)}{A(\varphi_2-\varphi_1)}} \ellipticF\left(\sqrt{\frac{v-\varphi_1}{\varphi_2-\varphi_1}},\frac{\varphi_2-\varphi_1}{\varphi_3-\varphi_1}\right),\text{ for } \varphi_1<v<\varphi_2,
		$$
		where $\ellipticF$ is the incomplete elliptic integral of the first kind \cite{handbookfunctions}, defined as
        %\notapan{OJO: He quitado el cuadrado de la $k$ en la definición, para que no haya número imaginario en las expresiones. No coincide con la notación de \cite{handbookfunctions}, pero la usa así el mathematica}
		$$
		\ellipticF(z,k)=\int_0^z \dfrac{ds}{\sqrt{(1-s^2) (1-k s^2)}}.
		$$

The corresponding solutions to equation \eqref{c6EqGenZK} are
\begin{equation}\label{solg3simplesa1}
	u(x,y,t)=\left(\varphi_2-\varphi_1\right)\jacobisn^2\left( \xi,\frac{\varphi_3-\varphi_1}{\varphi_2-\varphi_1}\right)+\varphi_1\jacobisn^4\left( \xi,\frac{\varphi_3-\varphi_1}{\varphi_2-\varphi_1}\right)
	\end{equation}
and
\begin{equation}\label{solg3simplesa2}
	u(x,y,t)=\varphi_1+\left(\varphi_2-\varphi_1\right)\jacobisn^2\left( \xi,\frac{\varphi_3-\varphi_1}{\varphi_2-\varphi_1}\right),
	\end{equation}
	where $\jacobisn(z;k)$ is the Jacobi elliptic sine function \cite{handbookfunctions}, and 
	$$
	\xi=\sqrt{\frac{A(\varphi_1-\varphi_2)}{12(M+N)}}(C_1-x-y+ct).
	$$
    Observe that in \eqref{solg3simplesa2} $c, C_1, \varphi_1,$ and $\varphi_2$ are free parameters (subject to the restrictions specified above), whereas $\varphi_3= \frac{3 c}{A}-\varphi_1-\varphi_2$.

	\item If $\frac{M+N}{A}>0$ then necessarily $\varphi_2<v<\varphi_3$ or $v<\varphi_1$. By performing, respectively, the change of variables $\tilde{v}=\sqrt{\frac{v-\varphi_2}{\varphi_3-\varphi_2}}$ or $\tilde{v}=\sqrt{\frac{\varphi_2-\varphi_1}{\varphi_2-v}}$, the following primitives are obtained
	$$
	H(v;\varphi_1,\varphi_2)=\sqrt{\frac{12(M+N)}{A(\varphi_2-\varphi_1)}} \ellipticF\left(\sqrt{\frac{v-\varphi_2}{\varphi_3-\varphi_2}},\frac{\varphi_2-\varphi_3}{\varphi_2-\varphi_1}\right),\text{ for } \varphi_2<v<\varphi_3,
	$$
	and
	$$
	H(v;\varphi_1,\varphi_2)=\sqrt{\frac{12(M+N)}{A(\varphi_2-\varphi_1)}} \ellipticF\left(\sqrt{\frac{\varphi_2-\varphi_1}{\varphi_2-v}},\frac{\varphi_2-\varphi_3}{\varphi_2-\varphi_1}\right),\text{ for } v<\varphi_1.
	$$

	The corresponding solutions to equation \eqref{c6EqGenZK} are
	\begin{equation}\label{solg3simplesb1}
		u(x,y,t)=\varphi_2+\left(\varphi_3-\varphi_2\right)\jacobisn^2\left( \xi,\frac{\varphi_2-\varphi_3}{\varphi_2-\varphi_1}\right)
	\end{equation}
	and
	\begin{equation}\label{solg3simplesb2}
		u(x,y,t)=\left(\varphi_1-\varphi_2\right)\jacobisn^2\left( \xi,\frac{\varphi_2-\varphi_3}{\varphi_2-\varphi_1}\right)+\varphi_2\jacobisn^4\left( \xi,\frac{\varphi_2-\varphi_3}{\varphi_2-\varphi_1}\right),
		\end{equation}
		where
		$$
	\xi=\sqrt{\frac{A(\varphi_2-\varphi_1)}{12(M+N)}}(C_1-x-y+ct).
	$$
Notice that $c, C_1, \varphi_1,$ and $\varphi_2$ remain as free parameters, and $\varphi_3= \frac{3 c}{A}-\varphi_1-\varphi_2$.

	\end{enumerate}

\end{itemize}

\subsection{\textbf{Case III}: \texorpdfstring{$B\neq 0$}{}}
%%%%%%%%%%%%%%%%%%%%%%%%%%%%%%%%%%%%%%%%%%%%%%5

We can rewrite the function $h(v;C_2,C_3)$ given in \eqref{h} as
\begin{equation}\label{hconlaK}
	h(v;C_2,C_3)=\sqrt{\dfrac{K}{v^4+\frac{2A}{B}v^3-\frac{6c}{B}v^2-\frac{C_2}{B}v+\frac{C_3}{B}}},
\end{equation}
where 
\begin{equation}\label{laK}
	K=\frac{-6 (M+N)}{B },
\end{equation}
and consider the following cases according to the multiplicities of the roots of the degree 4 polynomial \eqref{polinomioP}:

\begin{enumerate}
	\item  \textbf{$P(z)$ has a quadruple root}. This situation arises if and only if the wave speed is $c=-\frac{A^2}{4B}$, and the integration constants take the values
	$$
	C_2=-\frac{A^3}{2B^2}, \quad C_3=\frac{A^4}{16 B^3}.
	$$

	In this case, the quadruple root is $\varphi=-\frac{A}{2 B}$, and the function \eqref{hconlaK} becomes
	\begin{equation}\label{funintegrar1Quad}
		h(v)=\frac{\sqrt{K }}{\left(v+\frac{A}{2B}\right)^2}, \quad v\in \mathbb R,
	\end{equation}
	where necessarily $K>0.$ %in order to $h$ be a real-valued function. 
 A primitive of \eqref{funintegrar1Quad} is given by
	$$
H(v)=-\frac{\sqrt{K}}{v+\frac{A}{2B}}.
	$$
The corresponding solution to the modified ZK equation \eqref{c6EqGenZK}, obtained from \eqref{implicitTodas}, is 
	\begin{equation}\label{solg4cuadruple}
		u(x,y,t)=-\frac{A}{2 B}\pm\frac{\sqrt{K}}{C_1-x-y-\frac{A^2}{4B}t},
	\end{equation}
with $C_1$ a free parameter.

	\item  \textbf{$P(z)$ has a triple root}. It can be checked that this situation holds if and only if the wave speed $c$ satisfies $A^2+ 4 B c>0$ and the constants $C_2,C_3$ take the values
	$$
C_2=\frac{A(A^2+6 B c)\pm \sqrt{(A^2+4 B c)^3} }{B^2}
	$$
	and 
	$$
C_3=-\frac{A^2(A^2+6  B c)\pm A\sqrt{(A^2+ 4 B c)^3}+6 B^2 c^2}{2B^3}.
	$$

The roots are
\begin{equation}\label{raiceslatriple}
	\begin{aligned}
		\varphi_1&=- \frac{A\pm \sqrt{A^2+4 B c}}{2B} \text{ (triple)},		\\
		\varphi_2&=-\frac{A\mp 3\sqrt{A^2+4 B c}}{2B},
	\end{aligned}
\end{equation}
and it can be checked that $\varphi_1\neq \varphi_2$, since $c\neq-\frac{A^2}{4B}$.

In this case, the function \eqref{hconlaK} is
$$
h(v)= \sqrt{\frac{K}{(v-\varphi_1)^3(v-\varphi_2)} },
$$
where $v$ is such that $K(v-\varphi_1)(v-\varphi_2)>0$, and a primitive is given by
$$
H(v)=\dfrac{2}{\varphi_2-\varphi_1}\sqrt{K\frac{v-\varphi_2}{v-\varphi_1}}
$$
where $\frac{K(v-\varphi_2)}{v-\varphi_1}>0$.

The corresponding solution to equation \eqref{c6EqGenZK} is
	\begin{equation}\label{sol1RTriple}
		u(x,y,t)=\varphi_1+\frac{ \varphi_2 -\varphi_1}{1-\frac{(\varphi_2 -\varphi_1)^2}{4K}(C_1-x-y+ct)^2},
	\end{equation}
	where the constants $c$ and $C_1$ are the free parameters.
	
\item  \textbf{$P(z)$ has two different real double roots $\varphi_1$ and $\varphi_2$}. This occurs if and only if $c$ is such that $A^2+ 4 B c>0$ and the integration constants take the values
\begin{equation}\label{constantes2dobles}
    C_2=\frac{A\left(A^2+6 B c\right) }{B^2}, \quad C_3= \frac{\left(A^2+6 B c\right)^2}{4B^3}.
\end{equation}

The roots are
\begin{equation}\label{}
\begin{aligned}
	\varphi_1&=-\frac{ A+ \sqrt{3 (A^2+4 B c)}}{2B},\\	
	\varphi_2&=-\frac{ A- \sqrt{3 (A^2+4 B c)}}{2B},
\end{aligned}
\end{equation}
and the function \eqref{hconlaK} is 
\begin{equation}\label{funcionintegrar2dobles}
	h(v)= \sqrt{\frac{K}{(v-\varphi_1)^2(v-\varphi_2)^2} }, \quad v\in \mathbb R,
\end{equation}
where necessarily $K>0$.

The following primitive is obtained by means of the change of variables $\tilde{v}=\frac{v-\varphi_2}{v-\varphi_1}$
\begin{equation}\label{H2rD+}
H(v)=
\begin{cases}
	\dfrac{\sqrt{K}}{\varphi_1-\varphi_2}\ln \left( \dfrac{\varphi_2-v}{v-\varphi_1}\right),	& \text{if }\, \dfrac{v-\varphi_2}{v-\varphi_1} <0 \\
	\dfrac{\sqrt{K}}{\varphi_2-\varphi_1}\ln \left( \dfrac{v-\varphi_2}{v-\varphi_1}\right),	& \text{if }\, \dfrac{v-\varphi_2}{v-\varphi_1} >0 	
\end{cases}
\end{equation}

This allows us to obtain the following families of solutions for equation \eqref{c6EqGenZK}:
\begin{equation}\label{sol2rD+}
	u(x,y,t)=\varphi_1+\dfrac{\varphi_2-\varphi_1}{1+\exp\left({\mp \frac{\varphi_1-\varphi_2}{\sqrt{K}}(C_1-x-y+ct)}\right)},
\end{equation}
and
\begin{equation}\label{sol2rD-}
	u(x,y,t)=\varphi_1+\dfrac{\varphi_2-\varphi_1}{1-\exp\left({\mp \frac{\varphi_2-\varphi_1}{\sqrt{K}}(C_1-x-y+ct)}\right)}.
\end{equation}
In these expressions, $c$ and $C_1$ are the free parameters.

\item  \textbf{$P(z)$ has a pair of complex conjugate roots with double multiplicity, $a+b i$ and $a-bi$}.
This case arises when $c$ is such that $A^2+ 4 B c<0$, and the constants of integration take the same values as in \eqref{constantes2dobles}. The function to integrate \eqref{hconlaK} is 
\begin{equation}
	h(v)= \dfrac{\sqrt{K}}{(v-a)^2+b^2},
\end{equation}
where necessarily $K>0$, and $a$ and $b$ are, respectively the real and imaginary parts of the roots:
$$
a=-\frac{A}{2 B}, \quad b=\frac{\sqrt{-3(A^2+4Bc)}}{2B}.
$$

A primitive is
\begin{equation}
	H(v)=\dfrac{\sqrt{K}}{b} \arctan{\left(\frac{v-a}{b}\right)},
\end{equation}
and the solutions for equation \eqref{c6EqGenZK} are given by
\begin{equation}\label{solg4DoubleCompl}
	u(x,y,t)=a\mp b\tan\left(\frac{b}{\sqrt{K}}(C_1-x-y+ct)\right).
\end{equation}
We remark that $c$ and $C_1$ are free parameters.

\item  \textbf{$P(z)$ has a double root $\varphi_1$ and two other different real roots $\varphi_2$ and $\varphi_3$}. This possibility occurs if and only if the constants take the values
$$
C_2 = 2 \rho (2B \rho^2 + 3A \rho - 6c), \quad C_3 = \rho^2 (3B \rho^2 + 4A \rho - 6c)
$$
%for certain $\rho \neq -\frac{ A\pm \sqrt{3 (A^2+4 B c)}}{2B}$.
for certain $\rho \in \mathbb R$ satisfying $-2 B^2 \rho^2-2 A B \rho +A^2+6 Bc>0$.

The roots are
\begin{equation}\label{roots1RD}
	\begin{aligned}
		\varphi_1 &= \rho \, \text{(double)}, \\
		\varphi_2 &= -\frac{\rho B+A-\sqrt{-2 B^2 \rho^2-2 A B \rho +A^2+6 Bc} }{B}, \\
		\varphi_3 &= -\frac{\rho B+A+\sqrt{-2 B^2 \rho^2-2 A B \rho +A^2+6 Bc} }{B}.
	\end{aligned}
\end{equation}
Observe that $\varphi_2\neq \varphi_3$, because $-2 B^2 \rho^2-2 A B \rho +A^2+6 Bc\neq 0$.

In this case, the function \eqref{hconlaK} is
	\begin{equation}\label{funintegrar1RD}
		h(v;\rho)=\sqrt{\frac{K}{ (v-\varphi_1)^2(v-\varphi_2)(v-\varphi_3)} }, 
	\end{equation}
	where $v$ is such that $\frac{K}{ (v-\varphi_2)(v-\varphi_3)}>0$

	Without loss of generality, we may assume that $\varphi_2<\varphi_3$ (renaming the roots if necessary). Applying the change of variables
	$$
\tilde{v}=\sqrt{K \frac{v-\varphi_3}{v-\varphi_2}}
	$$

we obtain the following subcases:
\begin{enumerate}[label=(\alph*)]
	\item If $K (\varphi_1-\varphi_2)(\varphi_3-\varphi_1)>0$, a primitive is given by
	\begin{equation}
		H(v;\rho)=\dfrac{2K}{\sqrt{K(\varphi_1-\varphi_2)(\varphi_3-\varphi_1)}}\arctan\left(\sqrt{\dfrac{(\varphi_1-\varphi_2)(v-\varphi_3)}{(\varphi_3-\varphi_1)(v-\varphi_2)}}\right)
	\end{equation}
which yields the  solution to the mZK equation
	\begin{equation}\label{solg41D2realesA}
		u(x,y,t)=\varphi_2+\dfrac{\varphi_3-\varphi_2}{1-\frac{\varphi_3-\varphi_1}{\varphi_1-\varphi_2} \tan^2 \left(\frac{\sqrt{K (\varphi_3-\varphi_1)(\varphi_1-\varphi_2)}}{2K}(C_1-r)\right)},
	\end{equation}
	where $r=x+y-ct$.
	\item If $K (\varphi_1-\varphi_2)(\varphi_3-\varphi_1)<0$ then a primitive is given by
	\begin{equation}\nonumber
		\hspace{-1cm}
H(v;\rho) =
		\begin{cases} 
			\frac{2 K}{\sqrt{K (\varphi_2-\varphi_1) (\varphi_3-\varphi_1)}}  \arctanh\left(\sqrt{\frac{(\varphi_2-\varphi_1)(v-\varphi_3)}{(\varphi_3-\varphi_1)(v-\varphi_2)}}\right), & \text{if } \frac{(\varphi_2-\varphi_1)(v-\varphi_3)}{(\varphi_3-\varphi_1)(v-\varphi_2)}<1, \\
			\frac{-2 K}{\sqrt{K (\varphi_2-\varphi_1) (\varphi_3-\varphi_1)}}  \arctanh\left(\sqrt{\frac{(\varphi_3-\varphi_1)(v-\varphi_2)}{(\varphi_2-\varphi_1)(v-\varphi_3)}}\right), & \text{if } \frac{(\varphi_2-\varphi_1)(v-\varphi_3)}{(\varphi_3-\varphi_1)(v-\varphi_2)}>1.
		\end{cases}
		\end{equation}
\end{enumerate}

Thus we obtain two families of solutions to the mZK equation:
\begin{equation}\label{solg41D2realesB1}
	u(x,y,t)=\varphi_2+\dfrac{\varphi_3-\varphi_2}{1-\frac{\varphi_3-\varphi_1}{\varphi_2-\varphi_1} \tanh^2 \left(\frac{\sqrt{K (\varphi_3-\varphi_1)(\varphi_2-\varphi_1)}}{2K}(C_1-r)\right)},
\end{equation}
and
\begin{equation}\label{solg41D2realesB2}
	u(x,y,t)=\varphi_3+\dfrac{\varphi_2-\varphi_3}{1-\frac{\varphi_2-\varphi_1}{\varphi_3-\varphi_1} \tanh^2 \left(\frac{\sqrt{K (\varphi_3-\varphi_1)(\varphi_2-\varphi_1)}}{2K}(C_1-r)\right)},
\end{equation}
where $r=x+y-ct$.

The free parameters in \eqref{solg41D2realesA}, \eqref{solg41D2realesB1} and \eqref{solg41D2realesB2} are $c$, $C_1$, and $\rho$ because $\varphi_1, \varphi_2,\varphi_3$ are given by \eqref{roots1RD}.

\item  \textbf{$P(z)$ has a double root $\varphi_1$ and two complex conjugate roots $a\pm bi$}.

This case corresponds to the same values of the constants as in the previous case
$$
C_2=2 \rho (2B\rho^2+3A\rho-6c), \quad C_3= \rho^2(3 B\rho^2+4 A\rho-6c),
$$
but now the free parameter $\rho \in \mathbb R$ satisfies $-2 B^2 \rho^2-2 A B \rho +A^2+6 Bc<0$.

In this case, the roots are:
\begin{equation}\label{}
\begin{aligned}
	\varphi_1&=\rho \text{ (double)},\\
	\varphi_2&=a+ib,\\
	\varphi_3&=a-ib,	
\end{aligned}
\end{equation}
where the real and imaginary parts of the complex conjugate roots are
$$
a=-\frac{\rho B+A}{B},\quad b=-\frac{\sqrt{2 B^2 \rho^2+2 A B \rho -A^2-6 Bc} }{B},
$$
respectively.

The function to be integrated \eqref{hconlaK} is
	\begin{equation}
		h(v;\rho)=\sqrt{\frac{K}{ (v-\varphi_1)^2((v-a)^2+b^2)} },
	\end{equation}
	where necessarily  $K>0$.

	A primitive can be obtained by subsequently applying the changes of variables $v=a+b \sinh(\tilde{v})$, $\hat{v}=\tanh(\tilde{v}/2)$ and $\bar{v}=\frac{(\varphi_1-a) \hat{v}+b}{\sqrt{(\varphi_1-a)^2+b^2}}$:
	\begin{equation}\nonumber
		H(v;\rho) =
\begin{cases} 
    -\sqrt{\dfrac{4K}{(\varphi_1-a)^2+b^2}}\arctanh\left(\dfrac{\sqrt{(\varphi_1-a)^2+b^2}}{v-\varphi_1-\sqrt{(v-a)^2+b^2}}\right), & \text{if } v<\varphi_1, \\
    \sqrt{\dfrac{4K}{(\varphi_1-a)^2+b^2}}\arctanh\left(\dfrac{v-\varphi_1-\sqrt{(v-a)^2+b^2}}{\sqrt{(\varphi_1-a)^2+b^2}}\right),  & \text{if } v>\varphi_1.
\end{cases}
	\end{equation}

The corresponding solutions to \eqref{c6EqGenZK} are given by
\begin{equation}\label{solg4dobley2comp1}
	u(x,y,t)=\frac{a^2 + b^2 - \left( \varphi_1 \pm \sqrt{(\varphi_1-a)^2+b^2} \tanh\left(\xi\right)\right)^2}{2 \left(a-  \varphi_1 \mp \sqrt{(\varphi_1-a)^2+b^2} \tanh\left(\xi\right)\right)},
\end{equation}
and 
\begin{equation}\label{solg4dobley2comp2}
	u(x,y,t)=\pm\frac{  \left(\sqrt{(\varphi_1-a)^2+b^2} \pm \varphi_1 \tanh\left(\xi\right) \right)^2-(a^2 + b^2) \tanh^2\left(\xi\right)}{2 \tanh\left(\xi\right) \left(\sqrt{(\varphi_1-a)^2+b^2} \pm ( \varphi_1-a) \tanh\left(
\xi\right) \right)}
\end{equation}
where $\xi=\sqrt{\frac{ (\varphi_1-a)^2+b^2}{4K}}
(C_1 - x-y+ct)$. 

The free parameters in \eqref{solg4dobley2comp1} and \eqref{solg4dobley2comp2} are $c$, $C_1$, and $\rho$ because $a,b$ are given by \eqref{roots1RD}.

\item \textbf{$P(z)$ has four different roots $\varphi_1, \varphi_2,\varphi_3$ and $\varphi_4$}. This occurs if and only if the constants take the values
\begin{equation}\label{constantes4dist}
	\begin{aligned}
		C_2&= 4B\rho^3 + 6A\rho^2 + 4B\lambda\rho + 2A\lambda - 12c\rho,		\\
		C_3&=3B\rho^4 + 4A\rho^3 - 2B\lambda\rho^2 - 4A\lambda\rho - B\lambda^2 - 6c\rho^2 + 6c\lambda
	\end{aligned}
\end{equation}
for certain $\rho,\lambda \in \mathbb R$ satisfying:
$$
\lambda\neq 0,\quad A^2 + 6 B c-2 B^2 \rho^2 - 2 A B \rho - B^2 \lambda \neq 0.
$$

The roots are:
\begin{equation}\label{roots4dist}
	\begin{array}{l}
		\varphi_1=\rho+\sqrt{\lambda},\\
		\varphi_2=\rho-\sqrt{\lambda},\\
		\varphi_3=-\dfrac{B \rho + A}{B}+\sqrt{\dfrac{A^2 + 6 B c-2 B^2 \rho^2 - 2 A B \rho - B^2 \lambda }{B^2}},\\
		\varphi_4=-\dfrac{B \rho + A}{B}-\sqrt{\dfrac{A^2 + 6 B c-2 B^2 \rho^2 - 2 A B \rho - B^2 \lambda }{B^2}}.
	\end{array}
\end{equation}

Consider, first, the case when all the roots are real, that is, 
$$
\lambda> 0, \quad A^2 + 6 B c-2 B^2 \rho^2 - 2 A B \rho - B^2 \lambda> 0.
$$
The function \eqref{hconlaK} to be integrated is, in this case,
	$$
	h(v;\rho,\lambda)=\sqrt{\frac{K}{ (v-\varphi_1)(v-\varphi_2)(v-\varphi_3)(v-\varphi_4)} },
	$$
with $v$ such that $\frac{K}{ (v-\varphi_1)(v-\varphi_2)(v-\varphi_3)(v-\varphi_4)}>0$.

	To find a primitive, we rename the roots if necessary so that
	$$
\varphi_1<\varphi_2<\varphi_3<\varphi_4.
	$$
	%\notapan{Hay que distinguir 4 cambios distintos: $K>0$ y $K<0$, y dentro de cada uno hay que ajustar dominio para que ellipticF esté bien definida (argumento menor que 1)}
We consider the following subcases:
\begin{enumerate}[label=(\alph*)]
	\item $K>0$. In turn, this case is divided into two cases:
\begin{itemize}
		\item 	
 For $v<\varphi_1$ or $v>\varphi_4$, the change of variables
		$$
		\tilde{v}=\sqrt{\dfrac{(\varphi_4-\varphi_2)(v-\varphi_1)}{(\varphi_4-\varphi_1)(v-\varphi_2)}}
		$$
		allows us to obtain the primitive
   \begin{equation}%\label{masH}
\nonumber H(v;\rho,\lambda)=\sqrt{\frac{ 4 K}{(\varphi_3-\varphi_1)(\varphi_4-\varphi_2)}} \ellipticF\left(\sqrt{\frac{(\varphi_4-\varphi_2)\left(v-\varphi_1\right)}{(\varphi_4-\varphi_1)\left(v-\varphi_2\right)}} , \tilde{k} \right),
			\end{equation}
            where $\tilde{k}=\frac{(\varphi_3-\varphi_2)\left(\varphi_4-\varphi_1\right)}{(\varphi_3-\varphi_1)\left(\varphi_4-\varphi_2\right)}$.
        
			This primitive leads to the following family of solutions for the mZK equation:
			\begin{equation}\label{sol4dist1}
				u(x,y,t)=\varphi_2-\dfrac{\varphi_2-\varphi_1}{1-\dfrac{\varphi_4-\varphi_1}{\varphi_4-\varphi_2} \jacobisn^2\left( \xi , {\dfrac{(\varphi_3-\varphi_2)(\varphi_4-\varphi_1)}{(\varphi_3-\varphi_1)(\varphi_4-\varphi_2)}}\right)}
			\end{equation}
			where $\jacobisn(z;k)$ is the Jacobi elliptic sine function \cite{handbookfunctions}, and 
			$$
			\xi=\sqrt{\frac{(\varphi_3-\varphi_1)(\varphi_4-\varphi_2)}{4K}}(C_1-x-y+ct).
			$$	
    
		\item
        On the other hand, for $\varphi_2<v<\varphi_3$, we use the change of variables
		$$
		\tilde{v}=\sqrt{\dfrac{(\varphi_4-\varphi_1)(v-\varphi_2)}{(\varphi_4-\varphi_2)(v-\varphi_1)}}
		$$
		to obtain the primitive
        \begin{equation}%\label{masH2}
        \nonumber
H(v;\rho,\lambda)=\sqrt{\frac{4 K}{(\varphi_3 - \varphi_2)(\varphi_4 - \varphi_1)}} \, \ellipticF\left(  \sqrt{\frac{(\varphi_4 - \varphi_1)(v - \varphi_2)}{(\varphi_4 - \varphi_2)(v - \varphi_1)}}, k_2 \right)
		\end{equation} where
        $$k_2={\frac{(\varphi_3 - \varphi_1)(\varphi_4 - \varphi_2)}{(\varphi_3 - \varphi_2)(\varphi_4 - \varphi_1)}},$$
which yields to the following family of solutions:
		\begin{equation}\label{sol4dist2}
			u(x,y,t)=\varphi_1+\dfrac{\varphi_2-\varphi_1}{1-\dfrac{\varphi_4-\varphi_2}{\varphi_4-\varphi_1} \jacobisn^2\left( \xi , {\dfrac{(\varphi_3 - \varphi_1)(\varphi_4 - \varphi_2)}{(\varphi_3 - \varphi_2)(\varphi_4 - \varphi_1)}}\right)},
		\end{equation}
		where $\jacobisn(z;k)$ is the Jacobi elliptic sine function \cite{handbookfunctions}, and
		$$
		\xi=\sqrt{\frac{(\varphi_3-\varphi_2)(\varphi_4-\varphi_1)}{4K}}(C_1-x-y+ct).
		$$
	
\end{itemize}

\item $K<0$. We proceed as in the previous case to find the corresponding primitives, omitting the explicit expressions. We only indicate the changes of variables and the final solutions to the mZK equation, distinguishing the following cases: 

\begin{itemize}
	\item For $\varphi_1<v<\varphi_2$, the change of variables:
	$$
	\tilde{v}=\sqrt{\dfrac{(\varphi_4-\varphi_3)(v-\varphi_2)}{(\varphi_4-\varphi_2)(v-\varphi_3)}}
	$$
	leads to
	\begin{equation}\label{sol4dist3}
		u(x,y,t)=\varphi_3-\dfrac{\varphi_3 - \varphi_2}{1- \dfrac{\varphi_4 - \varphi_2}{\varphi_4 - \varphi_3}  \jacobisn^2\left( \xi, {\dfrac{(\varphi_3 - \varphi_1)(\varphi_4 - \varphi_2)}{(\varphi_2 - \varphi_1)(\varphi_4 - \varphi_3)}} \right)},
	\end{equation}
	where 
	$$
\xi=\sqrt{\frac{(\varphi_1-\varphi_2)(\varphi_4-\varphi_3)}{4K}}(C_1-x-y+ct).
	$$

	\item For $\varphi_3<v<\varphi_4$, the change of variables:
	$$
	\tilde{v}=\sqrt{\dfrac{(\varphi_4-\varphi_2)(v-\varphi_3)}{(\varphi_4-\varphi_3)(v-\varphi_2)}}
	$$
	leads to
	\begin{equation}\label{sold4dist4}
		u(x,y,t)=\varphi_2+\dfrac{\varphi_3 - \varphi_2}{1- \dfrac{\varphi_4 - \varphi_3}{\varphi_4 - \varphi_2}  \jacobisn^2\left(\xi, {\dfrac{(\varphi_2 - \varphi_1)(\varphi_4 - \varphi_3)}{(\varphi_3 - \varphi_1)(\varphi_4 - \varphi_2)}} \right)},
	\end{equation}
	where
	$$
\xi= \sqrt{\frac{(\varphi_1-\varphi_3)(\varphi_4-\varphi_2)}{4K}}(C_1-x-y+ct).
	$$	
 
\end{itemize}

\end{enumerate}

We remark that $c$, $C_1$, $\rho$, and $\lambda$ are the free parameters in expressions \eqref{sol4dist1}, \eqref{sol4dist2}, \eqref{sol4dist3}, and \eqref{sold4dist4}, since $\varphi_1,\varphi_2,\varphi_3,\varphi_4$ are given by \eqref{roots4dist}.

The case where
$$
\lambda < 0\quad 
\mbox{or}\quad
A^2 + 6Bc - 2B^2 \rho^2 - 2AB\rho - B^2 \lambda < 0,
$$
corresponds to complex roots. It can be checked that, by following the same approach as in the previous analysis (without reordering the roots), we obtain identical expressions. Therefore, these solutions are also valid for equation \eqref{c6EqGenZK}, even in the complex case.

\end{enumerate}

The families of solutions obtained in this section are classified in Table \ref{tablasoluciones}, where we indicate the number of the arbitrary  parameters in each family.

\renewcommand{\arraystretch}{1.5}%para aumentar la altura de las filas
\begin{table}
    \centering
    \begin{tabular}{|p{0.1\linewidth}|p{0.25\linewidth}|c|p{0.25\linewidth}|}
		
        \hline
        \textbf{Degree of $P$} & \textbf{Roots of $P$} & \textbf{Subcase} & \textbf{$k$-parameter family of solutions} \\
        \hline
        2 &  & $\frac{M+N}{c}>0$& $k=4$, \eqref{solg21} \\
        \cline{3-4}
        &&$\frac{M+N}{c}<0$& $k=4$, \eqref{solg22} \\
        \hline
        3& One triple root & &$k=2$,  \eqref{solg3triple}\\
        \cline{2-4}
           & One double root &$M+N>0$& $k=3$, \eqref{solg3caso2a1}, \eqref{solg3caso2a2} \\
		   \cline{3-4}
		   &         &$M+N<0$&$k=3$, \eqref{solg3caso2a1}, \eqref{solg3caso2b}\\
        \cline{2-4}
           & Three different roots &$\frac{M+N}{A}<0$&$k=4$, \eqref{solg3simplesa1}, \eqref{solg3simplesa2}\\
		   \cline{3-4}
		   & &$\frac{M+N}{A}>0$&$k=4$, \eqref{solg3simplesb1}, \eqref{solg3simplesb2}\\
        \hline
        4& One quadruple root &&$k=1$, \eqref{solg4cuadruple} \\
        \cline{2-4}
           & One triple root &&$k=2$, \eqref{sol1RTriple}\\
        \cline{2-4}
           & Two double real roots &&$k=2$, \eqref{sol2rD+}, \eqref{sol2rD-}\\
		   \cline{2-4}
		   & Two double complex roots &&$k=2$, \eqref{solg4DoubleCompl}\\
		   \cline{2-4}
		   & One double real root $\varphi_1$  and two single real roots $\varphi_2,\varphi_3$&$K (\varphi_1-\varphi_2)(\varphi_3-\varphi_1)>0$&$k=3$, \eqref{solg41D2realesA}\\
		   \cline{3-4}
		   & &$K (\varphi_1-\varphi_2)(\varphi_3-\varphi_1)<0$&$k=3$, \eqref{solg41D2realesB1}, \eqref{solg41D2realesB2}\\
		   \cline{2-4}
		   & One double real root  and two complex conjugate roots&&$k=3$, \eqref{solg4dobley2comp1}, \eqref{solg4dobley2comp2}\\
		   \cline{2-4}
		   & Four single roots& $K>0$&$k=4$, \eqref{sol4dist1}, \eqref{sol4dist2} \\
		   \cline{3-4}
		   & & $K<0$&$k=4$, \eqref{sol4dist3}, \eqref{sold4dist4}\\
        \hline
    \end{tabular}
    \caption{Classification of $k$-parameter families of traveling wave solutions of the mZK equation, based on the degree  of polynomial $P$ given in \eqref{polinomioP} and the multiplicity of its roots.}
    \label{tablasoluciones}
\end{table}

% \begin{center}
%     \begin{tabularx}{\textwidth}{|>{\hsize=0.3\hsize}X|>{\hsize=0.7\hsize}X|>{\hsize=1.5\hsize}X|>{\hsize=0.5\hsize}X|}
%         \hline
%         \textbf{Degree of $P$} & \textbf{Roots of $P$} & \textbf{Parameter restriction} & \textbf{Family of solutions} \\
%         \hline
%         2 & - & $\frac{M+N}{c}>0$& \eqref{solg21} \\
%         \cline{3-4}
%         &&$\frac{M+N}{c}<0$& \eqref{solg22} \\
%         \hline
%         3& Triple root & & \eqref{solg3triple}\\
%         \cline{2-4}
%            & Double root &$M+N>0$& \eqref{solg3caso2a1}, \eqref{solg3caso2a2} \\
%            \cline{3-4}
%            &         &$M+N<0$&\eqref{solg3caso2a1}, \eqref{solg3caso2b}\\
%         \cline{2-4}
%            & Single roots && (20), (21), (22), (23) \\
%         \hline
%         4& Quadruple root &&  \\
%         \cline{2-4}
%            & Double root && (16), (17), (18) \\
%         \cline{2-4}
%            & Single roots && (20), (21), (22), (23) \\
%         \hline
%     \end{tabularx}
% \end{center}

\section{Comparison with previous works}\label{section4}
%%%%%%%%%%%%%%%%%%%%%%%%%%%%%%%%%%%%%%
The modified Zakharov--Kuznetsov equation \eqref{c6EqGenZK} has been widely studied, yielding various particular solutions in the previous literature. This work unifies and generalizes results earlier obtained through  a diverse range  of powerful integration methods. 
Rather than providing an exhaustive identification of all specific solutions found so far, the analysis presented in this section aims to show how some of these solutions  corresponds to particular cases within the general families of solutions found in this work, arising from specific choices of the arbitrary parameters.

For example, in the recent paper \cite{mZKOsman} are found several traveling wave solutions to the modified Zakharov--Kuznetsov equation \eqref{c6EqGenZK} by using the double $(G'/G,1/G)$-expansion method. They are organized in three cases: 

\begin{itemize}
	\item In Case 1 of the paper \cite{mZKOsman}, the authors consider $c=-\lambda (M+N)$, with $\lambda>0$ a free parameter, and obtain the 1-parameter family of solutions
	\begin{equation}\label{osmancaso1}
		u(x,y,t)=-\frac{6 \lambda(M+N)}{\sqrt{A^2-6 B \lambda(M+N)} \cos (\sqrt{\lambda}(x+y-ct))+A},	
	\end{equation}
	given in equation (3.6) in \cite{mZKOsman}, which comes from the solution to \eqref{ZKode}:
	\begin{equation}\label{osmancaso1ODE}
		v(r)=-\frac{6 \lambda(M+N)}{\sqrt{A^2-6 B \lambda(M+N)} \cos (\sqrt{\lambda}r)+A}.	
	\end{equation}
	
To determine the family in Table \ref{tablasoluciones} that includes \eqref{osmancaso1}, we substitute \eqref{osmancaso1ODE} into equations \eqref{c6EqZKI3} and \eqref{c6EqZKI2}. This reveals that  this family corresponds to the values $C_2=C_3=0$. So in this case the polynomial \eqref{polinomioP} is
	$$
	P(z)=-Bz^4-2Az^3-6\lambda(M+N)z^2,
	$$
which has the following roots:
	  \begin{equation}
		\label{raicescaso1}
		\begin{split}
		\varphi_1 &= 0 \quad \text{(double)}, \\
		\varphi_2 &= \frac{-A - \sqrt{A^2 - 6 B \lambda (M+N)}}{B}, \\
		\varphi_3 &= \frac{-A + \sqrt{A^2 - 6 B \lambda (M+N)}}{B}.
		\end{split}
		\end{equation}

The authors only consider the case when %\notapan{Aunque en realidad no lo dicen, lo tenemos que deducir...} 
$A^2 - 6 B \lambda (M+N) >0$, which corresponds to  
one real double root and two other different real roots. Moreover, 
$$
K (\varphi_1-\varphi_2)(\varphi_3-\varphi_1)=\frac{36 \lambda (M+N)^2}{B^2}>0,
$$
so according to Table \ref{tablasoluciones}, equation \eqref{osmancaso1} must correspond to equation \eqref{solg41D2realesA}. Indeed, equation \eqref{osmancaso1} is obtained by substituting \eqref{laK} and \eqref{raicescaso1}, along with the particular value $C_1=0$, into the 3-parameter family of solutions \eqref{solg41D2realesA}.

On the other hand, still in the same section, the authors of \cite{mZKOsman} introduce in equation (3.7) the 1-parameter family of solutions
\begin{equation}\label{osmancaso1bis}
	u(x,y,t)=-\frac{6 \lambda(M+N)}{\sqrt{A^2-6 B \lambda(M+N)} \sin(\sqrt{\lambda}(x+y-ct))+A},	
\end{equation}
which fit again in the 3-parameter family of solutions \eqref{solg41D2realesA} with the substitutions \eqref{laK} and \eqref{raicescaso1}, but now with the particular value $C_1=\frac{\pi}{2\sqrt{\lambda}}$. 

	Note that  \cite{mZKOsman} does not consider  the case $A^2 - 6 B \lambda (M+N) < 0,$  which is however addressed in our analysis. In this case,  the polynomial has one real double root and a pair of complex conjugate roots, leading to the family of solutions \eqref{solg4dobley2comp1} and \eqref{solg4dobley2comp2}.

	\item In Case 2 in \cite{mZKOsman} it is assumed that $B(M+N)<0$. Also, it is considered that $\lambda=\frac{A^2}{6B(M+N)}$, which is equivalent to the assumption $c=-\frac{A^2}{6B}$, thereby obtaining the particular solutions (see  equation (3.10) in \cite{mZKOsman}):
	\begin{equation}\label{osmancaso2}
		\hspace{-0.5cm}u(x,y,t)=-\dfrac{A}{2B}\left(1\pm \coth(\sqrt{-\lambda}(x+y-ct))\pm \csch(\sqrt{-\lambda}(x+y-ct)) \right).
	\end{equation}
	
	%\notapan{Esta solución en realidad son 4, pues viene con $\pm$ en dos sitios. No sé si merece la pena identificar las cuatro constantes.}
	
	As before, by substituting the corresponding solution to the ODE \eqref{ZKode} into equations \eqref{c6EqZKI3} and \eqref{c6EqZKI2}, we deduce that they correspond to the choice $C_2=C_3=0$. Then, the polynomial \eqref{polinomioP} is
	$$
	P(z)=-Bz^4-2Az^3-\frac{A^2}{B}z^2,
	$$
	which has  two double real roots:
	  \begin{equation}
		\label{raicescaso2}
	\varphi_1 = 0, \quad \varphi_2 = -\frac{A}{B}.
		\end{equation}

According to Table \ref{tablasoluciones}, the four particular solutions  \eqref{osmancaso2} are special cases of  the 2-parameter families \eqref{sol2rD+} and \eqref{sol2rD-}. Indeed,  setting $C_1=0$ and substituting the expressions \eqref{laK} and \eqref{raicescaso2} into \eqref{sol2rD+} and \eqref{sol2rD-} yields the solutions \eqref{osmancaso2}.

% EN PARTICULAR, LAS SUSTITUCIONES DEBEN SER (VER MATHEMATICA 023):
% \begin{itemize}
% 	\item if we substitute \eqref{raicescaso2} into equation \eqref{sol2rD-}, together with $C_1=0$, we obtain the solution \eqref{osmancaso2} (case (+,+)). The sign in \eqref{sol2rD-} depends on the sign of $A$.
% 	\item if we substitute \eqref{raicescaso2} into equation \eqref{sol2rD+}, together with $C_1=0$, we obtain the solution \eqref{osmancaso2} (case (+,-)). The sign in \eqref{sol2rD-} depends on the sign of $A$.
% 	\item in \eqref{sol2rD+} for (-,+)
% 	\item in \eqref{sol2rD-} for (-,-)
% \end{itemize}

    \item In Case 3 of the same paper, the authors consider the value $c=-\dfrac{3 A^2}{16 B}$ and find a single solution (see equation (3.13) in \cite{mZKOsman}) given by
	\begin{equation}\label{solosmancaso3}
		u(x,y,t)=-\dfrac{3A}{4B}+\dfrac{24 A(M+N)}{A^2(x+y-ct)^2+24 B(M+N)}.
	\end{equation}

	By substituting the solution of \eqref{ZKode} corresponding to \eqref{solosmancaso3} into equations \eqref{c6EqZKI3} and \eqref{c6EqZKI2}, we deduce that this solution is framed in the choice $C_3=-\frac{27A^4}{256B^3}$ and $C_2=0$. Then, the polynomial \eqref{polinomioP} is
	$$
	P(z)=-Bz^4-2Az^3-\frac{9A^2}{8B}z^2+\frac{27A^4}{256B^3},
	$$
	with roots
	  \begin{equation}
		\label{raicescaso3}
		\varphi_1 = -\frac{3A}{4B}  \text{  (triple)},\quad \varphi_2 = \frac{A}{4B}.
		\end{equation}

		By inspecting Table \ref{tablasoluciones}, we deduce that this solution must be included in the 2-parameter family \eqref{sol1RTriple}, which, upon substitution of \eqref{laK} and \eqref{raiceslatriple}, yields
        $$
u(x,y,t)=\frac{-A + \sqrt{A^2 + 4Bc}}{2B} - \frac{12\sqrt{A^2 + 4Bc}(M + N)}{(A^2 + 4Bc)(x + y- ct-C_1)^2+6B(M + N) }.
        $$
It can be easily checked that \eqref{solosmancaso3} is recovered by taking $C_1=0$ and $c=-\dfrac{3 A^2}{16 B}$ in this 2-parameter family.
        
        %Indeed, setting $C_1=0$ and $c=-\dfrac{3 A^2}{16 B}$ and substituting \eqref{raicescaso3} and \eqref{laK} into equation \eqref{sol1RTriple} yields the particular solution \eqref{solosmancaso3}. 
        
\end{itemize}

  In summary, the 1-parameter families of solutions  \eqref{osmancaso1} and \eqref{osmancaso1bis} from \cite{mZKOsman} are special cases of  the 3-parameter families of solutions \eqref{solg41D2realesA}. Additionally, when  $A^2 - 6 B \lambda (M+N) < 0,$ the new solutions \eqref{solg4dobley2comp1} and \eqref{solg4dobley2comp2} arise. The four particular solutions  \eqref{osmancaso2} from \cite{mZKOsman}  are especial cases of the  2-parameter families \eqref{sol2rD+} and \eqref{sol2rD-}.  The single solution \eqref{solosmancaso3} is included  in the 2-parameter family \eqref{sol1RTriple}. Since the referred solutions \eqref{solg41D2realesA}, \eqref{solg4dobley2comp1}, \eqref{solg4dobley2comp2}, \eqref{sol2rD+},\eqref{sol2rD-} and \eqref{sol1RTriple} obtained in this work involve three  or two arbitrary parameters, it is clear that the number of new solutions has been significantly increased.

Another recent study \cite{alamin2022} investigates solitary wave solutions for the following variant of the mZK equation (see equation (7) in \cite{alamin2022}):
\begin{equation}\label{eqalamin}
	u_t+ u^2 u_x+ u_{xxx}+u_{xyy}=0,
\end{equation}
which corresponds to equation  \eqref{c6EqGenZK}  by setting $A=0$ and $B=M=N=1$. The auxiliary equation method, utilized in \cite{alamin2022}, produces a multitude of solutions that can be classified within our framework.

As an example, consider the 1-parameter family of solutions provided in equation (11) in \cite{alamin2022}:
\begin{equation}\label{solalamin11}
u(x,y,t)=\pm \sqrt{3}\left(\sqrt{m^2-4 l n}\tan\left(\frac{\sqrt{4 l n-m^2}}{2} (x+y-(4 l n-m^2)t)\right)\right),
\end{equation}
which is a traveling wave with $c=4 l n-m^2$. This expression corresponds to the following family of solutions to the ODE \eqref{ZKode}:
\begin{equation}\label{solalamin11ODE}
	v(r)=\pm \sqrt{3}\left(\sqrt{m^2-4 l n}\tan\left(\frac{\sqrt{4 l n-m^2}}{2} r\right)\right).
\end{equation}
If we substitute \eqref{solalamin11ODE} into expressions \eqref{c6EqZKI3} and \eqref{c6EqZKI2} we obtain the values
\begin{equation}\label{}
C_3=9\left(m^2-4 l n\right)^2,\quad C_2=0.
\end{equation}

Therefore, the polynomial \eqref{polinomioP} is, in this case,
$$
P(z)=z^4-6(4 l n-m^2)z^2+9\left(m^2-4 l n\right)^2,
$$
which has two double real roots:
\begin{equation}\label{alaminroots}
\varphi_1=-\sqrt{9(4 l n-m^2)},\quad  \varphi_2=\sqrt{9(4 l n-m^2)}.
\end{equation}

According to Table \ref{tablasoluciones}, the family \eqref{solalamin11} must be included in the 2-parameter families \eqref{sol2rD+} and \eqref{sol2rD-}. Indeed, substituting equations \eqref{laK} and \eqref{alaminroots} into \eqref{sol2rD+}, along with $c=4 l n-m^2$ and $C_1=0$, gives rise to equation \eqref{solalamin11}, as the reader can check.

A similar approach can be followed to classify the remaining solutions presented in \cite{alamin2022} for the particular mZK equation \eqref{eqalamin}.

\section{Some examples}\label{section5}
%%%%%%%%%%%%%%%%%%%%%%%%%%%%%%%%%%%%%%%%%%%%%%%%%%%%%%%%%%%%%%5

This section presents several  examples that demonstrate how the results of this study can be applied to identify interesting solutions for different variants of the mZK equation.

\subsection{Kink solutions}
%%%%%%%%%%%%%%%%%%%%%%%%%%%
% mathematica 023

The modified Zakharov--Kuznetsov equation 
\begin{equation}\label{eje2}
	u_t+ u u_x +2u^2 u_x- u_{xxx}-\frac{1}{3}u_{xyy}=0,
\end{equation}
corresponds to equation \eqref{c6EqGenZK} when $A=1,B=2,M=-1,N=-\frac{1}{3}$.
From equation \eqref{laK} we get  $K=-\frac{6(M+N)}{B}=4.$ 

If we take $C_2=\frac{1}{4}(1+12c)$ and $C_3=\frac{1}{32}(1+12c^2)$, the polynomial \eqref{polinomioP} has two double roots with values
\begin{equation}\label{}
\varphi_1= -\frac{1}{4}-\frac{1}{4}\sqrt{3+24c},\quad \varphi_2=-\frac{1}{4}+\frac{1}{4}\sqrt{3+24c}.
\end{equation}

For $c>-\frac{1}{8}$ both roots are real numbers so, according to Table \ref{tablasoluciones}, the solutions are given by  \eqref{sol2rD+} and \eqref{sol2rD-}. For example,  family of solutions \eqref{sol2rD+} with this choice of parameters reads as
\begin{equation}\label{soleje2}
	\begin{aligned}
		u(x,y,t)&=-\frac{1}{4}-\frac{\sqrt{3+24c}}{4} \tanh\left(\frac{\sqrt{3+24c}}{8} (C_1-x-y+ct)\right).
	\end{aligned}
\end{equation}
Different choices of the two parameters $C_1$ and $c>-1/8$ provide a large class of kink solutions for the mZK equation \eqref{eje2}. 
\subsection{Bright soliton solutions}
%%%%%%%%%%%%%%%%%%%%%%%%%%%%%%%%%%%%%%%%%%%
%mathematica 023

Consider the modified Zakharov--Kuznetsov equation 
\begin{equation}\label{}
	u_t+ u^2 u_x+ u_{xxx}+u_{xyy}=0,
\end{equation}
corresponding to equation \eqref{c6EqGenZK} for $A=0,B=M=N=1.$
If we set $C_2=C_3=0$, then the polynomial \eqref{polinomioP} becomes in this case
$$
P(z)=z^4-6cz^2
$$
with roots
\begin{equation}\label{}
\varphi_1= 0 \text{ (double)},\quad \varphi_2=\sqrt{6c},\quad \varphi_3=-\sqrt{6c}.
\end{equation}
From \eqref{laK}, we obtain $K=-\frac{6(M+N)}{B}=-12.$ Thus,
$$
K (\varphi_1-\varphi_2)(\varphi_3-\varphi_1)=-72 |c|<0.
$$
Assuming $c>0$ (so that $\varphi_2,\varphi_3$ are real), and according to Table \ref{tablasoluciones}, the traveling wave solutions for this choice of parameters are described by equations \eqref{solg41D2realesB1} and \eqref{solg41D2realesB2}. The substitution in these expressions of the corresponding values of the parameters yields
\begin{equation}\label{}
		u(x,y,t)=\mp\sqrt{6c}\sech\left(\sqrt{\frac{c}{2}}(C_1-x-y+ct)\right).
\end{equation}

This expression corresponds to bright soliton solutions (compare with equation (19) in \cite{ali2018lie}).

\subsection{Periodic solution}
%%%%%%%%%%%%%%%%%%%%%%%%%%

% mathematica 023

The modified Zakharov--Kuznetsov equation  \eqref{c6EqGenZK} corresponding to   the particular values $A=B=M=N=1$ becomes  
\begin{equation}\label{}
	u_t+ u u_x +u^2 u_x+ u_{xxx}+u_{xyy}=0,
\end{equation}
The value $K=-\frac{6(M+N)}{B}=-12$ is calculated  from equation \eqref{laK}.

In order to find periodic solutions of the form \eqref{sol4dist3}, we require $C_2,C_3$ to satisfy \eqref{constantes4dist}. For example, we can take $\rho=0$ and $\lambda=1$, and a wave speed $c=1$, so that, by \eqref{roots4dist},  the roots of the polynomial \eqref{polinomioP} are:
$$
\varphi_1=-1-\sqrt{6},\varphi_2=-1,\varphi_3=1,\varphi_4=-1+\sqrt{6}.
$$

Setting $C_1=0$ and substituting these values into \eqref{sol4dist3} yields the periodic solution
$$
u(x,y,t)=-1	+\frac{\sqrt{6}}{2 \jacobisn^2\left(\frac{x+y-t}{2\sqrt[4]{6}}\right)-1}.
$$

\section{Conclusions}
This article presents an exhaustive classification of all traveling wave solutions admitted by the modified Zakharov--Kuznetsov  equation, represented by equation \eqref{c6EqGenZK}. These solutions are derived from the general solution of the third-order ODE \eqref{ZKode}, using the recent \cinf-structure integration method. Integrating the final Pfaffian equation provided by the \cinf-structure requires finding a primitive  of the function \eqref{h}, expressed in terms of the polynomial $P$ given in \eqref{polinomioP}.

By analyzing the roots of $P$ and their multiplicities, we have classified all traveling wave solutions of the mZK equation into twenty-five distinct classes. Ten of these classes depend on four arbitrary parameters, nine involve three parameters, five depend on two, and only one class has one sole free parameter. Any traveling wave solution of the mZK equation necessarily belongs to one and only one of these twenty-five classes.

Additionally, the developed procedure allows us to identify the class to which any particular traveling wave solution obtained by other methods belongs (see Section \ref{section4}). This is done by first identifying the level set $\Sigma_{(C_2,C_3)}$ that defines the integral manifold of the second Pfaffian equation associated with the \cinf-structure. Once this is done, the polynomial $P$ can be calculated, and by examining its roots and their multiplicities, we can use Table \ref{tablasoluciones} to identify the family to which the given solution belongs. Furthermore, due to the presence of arbitrary parameters, many other particular solutions can be obtained by assigning specific values to these constants.

Finally, classifying the solutions into $k$-parameter families facilitates the search for relevant types of traveling wave solutions for any mZK equation of the form \eqref{c6EqGenZK}, such as bright solitons, kink solutions, periodic solutions, etc. This is especially useful in improving our understanding of the many physical phenomena that can be modeled by different versions of the mZK equation.

\section*{Acknowledgments}

 The authors appreciate the partial financial support by the grant ``Operator Theory: an interdisciplinary approach", reference ProyExcel$\_$00780, a project financed in the 2021 call for Grants for Excellence Projects, under a competitive bidding regime, aimed at entities qualified as Agents of the Andalusian Knowledge System, in the scope of the Andalusian Research, Development and Innovation Plan (PAIDI 2020). Counseling of University, Research and Innovation of the Junta de Andaluc\'ia. 
 
 The authors also acknowledge the financial support of the Research Project PR2023-024 from {\it   Plan Propio de Est\'imulo y Apoyo a la Investigaci\'on y Transferencia 2022-2023} of the University of C\'adiz; as well as  the support from {\it   Junta de Andaluc\'ia} (Spain)  to the  research group FQM--377.

\bibliographystyle{unsrt}  
\bibliography{references.bib}

\end{document}